\documentclass[12pt]{amsart}

\usepackage{amssymb}

\usepackage{enumerate}

\usepackage{graphicx}
\usepackage{color}
\makeatletter
\@namedef{subjclassname@2010}{%
  \textup{2010} Mathematics Subject Classification}
\makeatother

\newtheorem{thm}{Theorem}[section]
\newtheorem{cor}[thm]{Corollary}
\newtheorem{lem}[thm]{Lemma}
\newtheorem{prop}[thm]{Proposition}

\newtheorem{mainthm}[thm]{Main Theorem}

\theoremstyle{definition}
\newtheorem{defin}[thm]{Definition}
\newtheorem{rem}[thm]{Remark}
\newtheorem{exa}[thm]{Example}

\numberwithin{equation}{section}

\begin{document}

%%%%% To ease editing, for IMPAN journals add:

\baselineskip=17pt

%%%%%%%%%%%

%% In the running head, replace first names by initials
%% and give an abbreviation of the title.

\title[On The  Rank Of  Congruent Elliptic Curves ]{On The  Rank Of  Congruent Elliptic Curves}

\author[F. Izadi]{Farzali Izadi}
\address{Farzali Izadi \ Department of Mathematics, \ Faculty of Science, \ Urmia University, \ Urmia 165-57153, Iran}
\email{f.izadi@urmia.ac.ir}

\author[H. R. Abdolmaleki]{Hamid Reza Abdolmaleki}
\address{H. R. Abdolmaleki, \  Department of  Mathematics, \  Faculty of Science,\  Azarbaijan Shahid Madani University,\\Tabriz 53751-71379, Iran}
\email{Hamid.Abdolmalki@gmail.com}

\date{}

\begin{abstract}
In this paper, $p$ and $q$ are two different odd primes. First, We construct the  congruent elliptic curves corresponding to  $p$, $2p$, $pq$, and $2pq,$ then, in the cases of congruent numbers, we determine the rank of the corresponding congruent elliptic curves.
\end{abstract}


\maketitle
\section{introduction}
The rank of an elliptic curve is a measure of the size of the set of rational points. However, the question is to ask how one can compute the exact size of this set of rational points. On the other hand, it is easy to find the rational points of a projective line or a plane curve defined by a quadratic equation.
Having said that, there is no known guaranteed algorithm to determine the rank and it is not known which numbers can occur as the rank of an elliptic curve. (see \cite{nag}).

The rational number $n$ is a congruent number if there are positive rational numbers $a, b, c$ such that $a^2+b^2 = c^2$ and $\frac{1}{2}ab= n$ , equivalently, there is a Pythagorean triangle with rational sides and the area equals to $n$.

In a modern language, $n$ is a congruent number if and only if the elliptic curve $ E: y^2= x^3 - n^2 x $ contains a rational point with $y\neq0$, equivalently, a rational point of infinite order, i.e., the rank of $E$ which is denoted by $r(E)$ is nonzero, (see {\cite{was}).

First we quote from Monsky \cite {mon}. Let $p_{1}, p_{3}, p_{5}$ and $p_{7}$ denote primes $\equiv 1, 3, 5$ and $7 \pmod 8$. Heegner \cite{hee}, and Brich \cite{B1}, proved that $2p_{3}$ and $2p_{7}$ are congruent numbers. Heegner asserted without proof that $p_{5}$ and $p_{7}$ are congruent numbers, and this claim is repeated in \cite{ste}. Monsky \cite{mon} has given a unified proof that the following are all congruent numbers :

$(1) : p_{5},\ p_{7},\ 2p_{7}$, and $2p_{3}$,

$(2) : p_{3} p_{5},\ p_{3} p_{7},\ 2p_{3} p_{5}$ and $2p_{5}p_{7}$,

$(3) : p_{1} p_{5}$ provided $(\frac{p_{1}}{p_{5}})=-1$, and $2p_{1} p_{3}$ provided
$(\frac{p_{1}}{p_{3}})=-1$ and $p_{1} p_{7}$, $2p_{1} p_{7}$, provided $(\frac{p_{1}}{p_{7}})=-1$,

In other words, any $N= 5, 6$ or $7\pmod8$ having at most 2 odd prime factors is a congruent number, with the following possible exception $N=pq$ or $2pq$ with $p\equiv1\pmod8$ and $(\frac{p_{1}}{p_{3}})=+1$.

\begin{rem} For $N$ in the form of $(1), (2)$ and $(3)$, our proof shows that the rank of $E^{\mathbb{Q}} _{N} : Y^2=X^3-N^2X$ is $1$. However, for $N=pq$ or $2pq$ with $p\equiv1\pmod8$ and $(\frac{p}{q})=+1$, this result is not true. For example, for $N=(521).(5)$ , the rank of $ E^{\mathbb{Q}} _{N}$ is 3.
\end{rem}
In section two, we recall the 2-descent method which is a classical method  for finding the rank of elliptic curves. Section three includes our results and contains five parts.\\
 In part (1), Let $E: y^2= x^3 - p^2 x$, be the corresponding congruent number elliptic curve for $p$. It is known that if $p\equiv 5, 7 \pmod 8$ then $p$ is congruent number (see\cite {mon}). Using the 2-descent method, we show in this two cases, $r(E)=1$. Moreover if $p\equiv 3 \pmod 8$ then $r(E)= 0$, i.e., $p$ is not congruent.

In part (2), for $p\equiv1\pmod 8$, we investigate that $r(E)= 2$. We show this is happen whenever $p$ satisfy in the two following conditions:

(1) $\exists a , b \in \mathbb{N}$: $p= a+b~, ~ a - b= \square$ and $a^2+b^2= \square$.

(2) $\exists a , b \in \mathbb{N}$: $p= a^2+b^2~, ~ (a, 2b)=1$ and $(a \pm 2b)^2+b^2= \square$.

In part (3), Let $E: y^2= x^3 -4 p^2 x$, be the corresponding congruent number elliptic curve for $2p$. It is known that if
$ p\equiv 3, 7 \pmod 8$ then $2p$ is congruent number (see\cite {mon}). Using the $2$-descent method, we show in this two case, $r(E)=1$. Moreover if $ p\equiv 5 \pmod 8$ then $r(E)= 0$, i.e., $2p$ is not congruent.

In part (4), Let $ E: y^2= x^3 - p^2q^2 x$, be the corresponding congruent number elliptic curve for $pq$. It is known that if
$[p\equiv 3$ and $q\equiv 5, 7\pmod 8]$ or $[p\equiv1$ and $ q\equiv 5, 7 \pmod 8$ such that $(\frac{p}{q})=-1]$ then $pq$ is a congruent number (see\cite {mon}). Using 2-descent method, we show in this four cases, $r(E)=1$. Moreover if $p, q\equiv3 \pmod 8$ then $r(E)= 0$, i.e., $pq$ is not  congruent.

In part (5), Let $E : y^2 = x^3 - 4p^2 q^2 x$. We know If $[p\equiv 5$ and
$q \equiv 3, 7\pmod8]$ or $[ p\equiv 1$ and $q \equiv 3, 7\pmod8$ such that $( \frac{p}{q}) = -1]$ then $2pq$ is congruent number (see\cite {mon}). Using the $2$- descent method, we show in this four cases, $r(E) =1$. Moreover if $[p, q\equiv 5\pmod8]$ or $[p \equiv1$ and $q \equiv5 \pmod8$ such that $(\frac{p}{q}) = -1],$ then $r(E) = 0$, i.e., $2pq$ is not congruent.

$~$
\section{2-descent method}
In this section we describe the $2$- descent method for determining the rank of an elliptic curve, (see \cite{coh}, Chapter 8 for more details). Let $E(\mathbb{Q})$ is the group of rational points on the elliptic curve $ E: y^2= x^3 + ax^2+bx$, where $a, b\in\mathbb{Q}$. Let $\mathbb{Q^*}$ is the multiplicative group of nonzero rational numbers and $\mathbb{Q^*}^2$ is the subgroup of squares of elements of $\mathbb{Q^*}$. Define the $2$-descent homomorphism $\alpha$ from $E(\mathbb{Q})$ to $\frac{\mathbb{Q^*}}{\mathbb{Q^*}^2}$ as follows :
\begin{equation*}
\alpha(P)=\left \{\begin{array}{ll}
1 ~ \pmod{\mathbb{Q^*}^2} & \text{ if}\ ~ P=\mathcal{O}= \infty, \\
b ~ \pmod{\mathbb{Q^*}^2} & \text{ if}\ P= (0,0), \\
x ~ \pmod{\mathbb{Q^*}^2 }& \text{ if}\ P= (x,y), x \neq 0. \end{array}\right.
\end{equation*}
Similarly, take the isogenous curve $\overline {E}: y^2= x^3 -2ax^2+(a^2-4b)x$ with the group of rational points $\overline{E}(\mathbb{Q})$. The 2-descent homomorphism $ \overline{\alpha} $ from $\overline{E}(\mathbb{Q})$ to $\frac{\mathbb{Q^*}}{\mathbb{Q^*}^2}$ as follows :
\begin{equation*}
\overline{\alpha}(P)=\left \{\begin{array}{ll}
1 \pmod{ \mathbb{Q^*}^2}& \text{ if}\ \overline{P}= \mathcal{O}= \infty, \\
a^2-4b \pmod{\mathbb{Q^*}^2} & \text{ if}\ \overline{P}= (0,0), \\
x \pmod {\mathbb{Q^*}^2} & \text{ if}\ \overline{P}= (x,y), x \neq 0. \end{array}\right.
\end{equation*}
The rank of $E(\mathbb{Q})$ which is denoted by $r$ is determined by
\begin {equation}\label{eq2}
2^r=\frac{\mid Im\alpha\mid\cdot \mid Im\overline{\alpha}\mid}{4}.
\end{equation}
The group $\alpha (E(\mathbb{Q}))$ equals to the classes modulo squares of $1,\ b$ and the positive and negative divisors of $b$ such that
\begin {equation*}
N^2= b_{1}m^4+am^2e^2+\frac{b}{b_{1}}e^4,
\end{equation*}
\begin {equation*}
\gcd(m, e)=\gcd(m, N)=\gcd(e, N)=\gcd(b_{1}, m)=\gcd(\frac{b}{b_{1}}, e)=1, me\neq 0.
\end{equation*}
If $(m, N, e)$ is a solution, then $ P=(\dfrac{b_{1}m^2}{e^2},\dfrac{b_{1}mN}{e^3})$ belongs to $E(\mathbb{Q})$ and
we have the same for $\overline{\alpha}$ as well.

$~$
\section{our results}

\subsection*{Part (1)}
$~$

According the 2-descent method, for the elliptic curves
\begin{equation}\label{eqA}
E: y^2=x^3 - p^2 x
\end{equation}
\ \ \ \ and
\begin{equation}\label{eqB}
\overline{E }: y^2=x^3 +4 p^2q^2 x.
\end{equation}
We have respectively
${\{\pm1}\} \subseteq Im\alpha \subseteq {\{\pm 1,\pm p}\}$ and
$ {\{1}\} \subseteq Im \overline{\alpha} \subseteq {\{1, 2, p, 2p}\}$.
Therefore, according to \eqref{eq2}, the maximum rank of \eqref{eqA} is 2. Moreover, the homogeneous equation of $E$ is

$(a_{1}): N^2= \pm(pm^4 -pe^4),\quad  \gcd(m, pe)=\gcd(e, pm)=\gcd(N, me)=1$,

and the homogeneous equations of $\overline{E}$ are
\begin{equation*}
\begin{array}{l}
(b_{1}): N^2= pm^4 +4pe^4,\quad  \gcd(m, 2pe)=\gcd(e, pm)=\gcd(N, me)=1,\\
\\
 (b_{2}): N^2= 2m^4 +2p^2e^4,\quad \gcd(m, 2pe)=\gcd(e, 2m)=\gcd(N, me)=1,\\
\\
 (b_{3}): N^2=2 pm^4 +2pe^4,\quad \gcd(m, 2pe)=\gcd(e, 2pm)=\gcd(N, me)=1.
\end{array}
\end{equation*}
First we study the solvability of the above homogeneous equations.

\begin{prop}\label{p31}
If $(b_{1})$ has integer solution then $p\equiv1 \pmod 4$.
\end{prop}
\begin{proof}
If $(b_{1})$ has integer solution then $N^2= p(m^4 +4e^4)$, $\gcd(m, p)=1$. There is an integer $u$ such that $pu^2=m^4+4e^4$. Hence $m^4 \equiv -4e^4\pmod p$.
 Therefore $(2e^2 {m^*}^2)^2 \equiv -1 \pmod p$. So $(\frac{-1}{p})=+1$, i.e., $p\equiv 1 \pmod 4$.
\end{proof}

\begin{prop}\label{p32}
If $(b_{2})$ has integer solution then $p\equiv\pm 1 \pmod 8$.
\end{prop}
\begin{proof}
If $(b_{2})$ has integer solution then $N^2= 2(m^4 +p^2e^4)$, $\gcd(m, p)=1$.
There is an integer $u$ such that $2u^2=m^4+p^2e^4$. So $2u^2\equiv m^4\pmod p$. Hence $(2um^{*^2})^2 \equiv 2 \pmod p$. So $(\frac{2}{p})=+1$, i.e., $p\equiv\pm 1 \pmod 8$.
\end{proof}

\begin{rem}\label{r33}
If $u$ is an integer number then $u^2 \equiv 0, 1, 4 \pmod 8$. Hence, $\gcd(u, 8)=1$ if and only if $u^2 \equiv 1 \pmod 8$.
\end{rem}

\begin{prop}\label{p34}
If $(b_{3})$ has integer solution then $p\equiv1\pmod 8$.
\end{prop}
\begin{proof}
 If  $(b_{3})$  has integer solution then $N^2= 2p(m^4  + e^4)$, $\gcd (em , 2)=1$.
There is an integer $u$ such that $2pu^2 = m^4 + e^4$. Hence $2pu^2 \equiv 2 \pmod{16}$. So $pu^2 \equiv 1\pmod8 $. We have $(u, 8)=1$. Therefore $u^2 \equiv 1 \pmod 8$. Consequently $p \equiv 1 \pmod 8$.
 \end{proof}

\begin{cor}\label{c35}
The above propositions are summarized in the followings:
\begin{enumerate}
\item $\pm p \in Im\alpha$ if and only $(a_{1})$ has integer solution.

\item $p \in Im\overline{\alpha}$ if and only if $ (b_{1})$ has integer solution. In this case $p\equiv 1, 5\pmod 8$.

\item $2 \in Im\overline{\alpha}$ if and only if $(b_{2})$ has integer solution.  In this case $p\equiv 1, 7 \pmod 8$.

\item $2p \in Im\overline{\alpha}$ if and only if $ (b_{3})$ has integer solution. In this case $p \equiv 1\pmod 8$.

\end{enumerate}
\end{cor}

\begin{cor}\label{c36}
By using previous corollary, we have:
\begin{enumerate}
\item Let $p\equiv 3\pmod 8$.

Then $(b_{1}),(b_{2}),(b_{3})$ do not have  any integer solutions, so $Im\overline{\alpha}= \lbrace1\rbrace$.

As \eqref{eq2}, implies $\mid Im\alpha\mid\geq4$.

So $Im{\alpha}= \lbrace\pm 1, \pm p\rbrace$. Hence $r = 0$, i.e., $p$ is not a congruent.
\\
\item Let $p\equiv 5 \pmod 8$.

It is possible for $(b_{1})$ to have integer solution. So $Im\overline{\alpha}\subseteq \lbrace1, p \rbrace$.

$p$ is a congruent number, (see \cite{mon}), we have $r\geq1$.

As \eqref{eq2}, implies $Im{\alpha}= \lbrace\pm 1, \pm p\rbrace$ and $Im\overline{\alpha} = \lbrace 1, p\rbrace$. Hence $r = 1$.
\\
\item Let $p\equiv 7 \pmod 8$.

It is possible for $(b_{2})$ to have integer solution. So $Im\overline{\alpha}\subseteq \lbrace1, 2 \rbrace$.

$p$ is a congruent number, (see \cite{mon}), we have $ r\geq1$.

 As \eqref{eq2}, $Im{\alpha}= \lbrace\pm 1, \pm p\rbrace$ and $Im\overline{\alpha} = \lbrace 1, 2\rbrace$.
Hence $r = 1$.
\end{enumerate}
\end{cor}
$~$

\subsection*{Part (2)}
$~$

In this section, for $p\equiv1\pmod 8$ we investigate that $r(E)=  2$.

$\bullet$ First we study the solvability of the homogeneous equation $(b_{1})$. We need some definitions.

\begin{defin}
$p$ is a $\alpha _{-}$Pythagorean whenever, there are coprime integer numbers $a, b$ and $c$ such that
$pc^2=a^2+b^2$.
\end{defin}

\begin{defin}
$p$ is a $\alpha_{-}^{-}$Pythagorean whenever there are coprime integers $a, b$ and $c$ such that  $pc^2=a^2+b^2$ and $(a - 2b)^2+b^2= \square$.
\end{defin}
\begin{defin}
$p$ is a $\alpha _{-}^{+}$Pythagorean whenever there are coprime integers $a, b$ and $c$ such that  $pc^2=a^2+b^2$ and $(a+2b)^2+b^2= \square$.
\end{defin}

\begin{defin}
$p$ is a $\alpha _{-}^{\pm}$Pythagorean whenever $p$ is a
$\alpha_{-}^{-}$Pythagorean or $\alpha_{-}^{+}$Pythagorean.
\end{defin}

\begin{rem}
Considering the results of part 2, we know $p\equiv 5 \pmod 8$ is $\alpha_{-}^{\pm}$Pythagorean.
\end{rem}

\begin{exa}
$37$ with $(a, b, c) = (22, 21, 5)$ is a $\alpha_{-}^{-}$Pythagorean.
\end{exa}

\begin{exa}
$41$ with $(a, b, c) = (5, 4, 1)$ is a $\alpha_{-}^{-}$Pythagorean.
\end{exa}

\begin{exa}
$149$ with $(a, b, c)= (10, 7, 1)$ is a $\alpha_{-}^{+}$Pythagorean.
\end{exa}

\begin{prop}
$(b_{1})$ has integer solution if and only if $p$ is a $\alpha_{-}^{\pm}$Pythagorean.
\begin{proof}
If $(b_{1})$ has integer solution then $N^2= p (m^4 + 4e^4), gcd(m, 2e)=1$. There is an integer number $u$ such that $pu^2= m^4+4e^4$.

$pu^2=(m^2 - 2me+2e^2) (m^2+2me+2e^2)$.

As $gcd(m^2 - 2me+2e^2, m^2+2me+2e^2)=1$, there are coprime integers $c$ and $w$ such that
\begin{enumerate}
\item $m^2-2me+2e^2= pc^2$ and $m^2+2me+2e^2= w^2$. So $(m-e)^2+e^2= pc^2$ and
 $[(m-e)+2e]^2+e^2= w^2$. Hence $p$ is a $\alpha_{-}^{+}$Pythagorean.

\item $m^2-2me+2e^2= w^2$ and $m^2+2me+2e^2= pc^2$.
So $(m+e)^2+e^2= pc^2$
\end{enumerate}
\qquad \quad and $[(m+e) - 2e]^2+e^2= w^2$. Hence $p$ is a $\alpha_{-}^{-}$Pythagorean.
\end{proof}
\end{prop}

\begin{rem}
If $p\equiv 1\pmod4\equiv1 , 5\pmod 8$ then there are unique positive integers $a$ and $b$ such that $p=a^2+b^2$.

$p \equiv 1 \pmod 8$ if and only if there are unique integers $k$ and $t$ such that $p= 16k^2+t^2, (t, 2k)=1$.

$p \equiv 5 \pmod 8$ if and only if there are unique integers $k$ and $t$ such that $p= 4k^2+t^2, (t , 2k)=(k , 2)=1$.
\end{rem}

\begin{lem}
If $p$ with $(a, b, c)$ is a $\alpha_{-}^{\pm}$Pythagorean, then

$a^2 \equiv 1$ ~,~ $b^2 \equiv 0 \pmod 8  \qquad\text {or}\qquad a^2 \equiv 4$  ~,~ $b^2 \equiv 1 \pmod 8$.
\begin{proof}
We have $pc^2=a^2+b^2 , (a , b)=1$ and $(a\pm2b)^2+b^2= \square$.

Since $pc^2$ is sum of two primitive squares, then $p$ and all factors of $c$ are form $4k+1$, where $k$ is
integer, (see\cite{Adl}). We know squares in mod 8 are $0, 1$ or $4$.

Suppose otherwise if $a^2 \equiv 0 , b^2\equiv 1$ or $a^2\equiv 1 , b^2\equiv 4 \pmod 8$,
then $\square= a^2+5b^2 \pm 4ab\equiv 5 \pmod 8$, which  is a contradiction.
\end{proof}
\end{lem}

\begin{cor}
Suppose $p$ is a $\alpha_{-}^{\pm}$Pythagorean:
\begin{enumerate}
\item If $p=4k^2+t^2\equiv 5\pmod 8$, then $a=2k$ and $b= t$, where $kt$ is odd.

\item If $p=16k^2+t^2\equiv 1\pmod 8$, then $a= t$ and $b= 4k$, where $t$ is odd.
\end{enumerate}
\end{cor}
$\bullet$  Next theorem tell us that the prime number $p \equiv 1 \pmod 8$ is a $\alpha_{-}^{\pm}$Pythagorean if and only if $c=1$, (see the above Definition).
\begin{thm}
(1) If $p= a^2+b^2$ and $2\parallel b$, then $p$ is congruent number.

(2) If $4 \mid b$ and $(a\pm2b)^2+b^2\neq \square$, then $p$ is not a $\alpha_{-}^{\pm}$Pythagorean.

\begin{proof}

(1) If $2 \parallel b$, then there is odd integer number $b_{0}$ such that $b = 2b_{0}$. So $p=a^2+4b_{0}^2\equiv 5 \pmod 8$, therefore $p$ is congruent number, (see\cite{mon}).

(2) If $4 \mid b$, $(a\pm2b)^2+b^2\neq \square$ and Suppose otherwise $p$ is a $\alpha_{-}^{\pm}$Pythagorean then there are coprime integers $a_{0},~ b_{0}$ and $c_{0}\neq 1$ such that
\begin{center}
  $ pc_{0}^2=a_{0}^2+b_{0}^2$\qquad\text{and}\qquad $4a_{0}^2+5b_{0}^2 \pm4a_{0}b_{0}=\square$.
\end{center}
 Therefore $pc_{0}^2$ is sum of two primitive numbers, then $p$ and all factors of $c_{0}^2$ are form $4k+1$. Consequently there are integers $m$ and $n$ such that $m$ is odd number and $n$ is nonzero even number such that $c_{0}^2=m^2+n^2$, then
\begin{center}
$pc_{0}^2= (b^2+a^2)(m^2+n^2)= (bm+an)^2+(bn-am)^2$.
\end{center}
We have $[a_{0}= bm+an ~,~ b_{0}= bn-am]$  or $[a_{0}= bm-an$~,~$ b_{0}= bn+am]$.

If $a_{0}= bm+an$ and $b_{0}= bn-am$, then $a^2n^2 \pm5a^2m^2\equiv\square \pmod 8$.

If $n^2\equiv 4$ and $m^2\equiv 1\pmod 8$, then $c_{0}^2\equiv5 \pmod 8$, which is a contradiction.

If $n^2\equiv 0, m^2\equiv 1$ and $a^2\equiv 1\pmod 8$, then $\square\equiv \pm 5 \pmod 8$, which is a contradiction.

proof for $a_{0}= bm-an$ and $b_{0}= bn+am$ is similar.
\end{proof}
\end{thm}
\begin{exa}
$17=16\times 1^2+1^2\equiv 0 + 1\pmod 8$, then $a^2=1$ and $b^2=16$. As $(a\pm 2b)^2+b^2\neq \square$, hence $17$ is not a $\alpha_{-}^{\pm}$Pythagorean.
\end{exa}
\begin{exa}
$41=16\times 1^2+5^2$, with $(a, b, c)= (5, 4, 1)$, is a $\alpha_{-}^{-}$Pythagorean.
\end{exa}

$\bullet$ Now, we study the solvability of the homogeneous equation $(b_{2})$, however we first need a definitions.

\begin{defin}
$p$ is a $\beta _{-}$Pythagorean whenever there are integers $a, e, m$ and $u$ such that $pe^2m^2= 2a^2 - u^2$ , $pe^2= 2a - m^2$ and $(e , m)=1$.
\end{defin}

\begin{rem}
Considering the results of part 2, we know $p\equiv 7 \pmod 8$ is $\beta_{-}$Pythagorean.
\end{rem}

\begin{prop}
$(b_{2})$ has integer solution if and only if $p$ is a $\beta _{-}$Pythagorean.
\begin{proof}
If $(b_{2})$ has integer solution then $N^2= 2(m^4+p^2e^4), gcd(m, e)=1$.

$N^2= 2 (m^4+p^2e^4)$ if and only if there is a integer number $u$ such that $2u^2=m^4+p^2e^4$
if and only if $(m^2+pe^2)^2= 2(u^2+pm^2e^2)$ if and only if there is an integer number $a$ such that $u^2+pm^2e^2 = 2a^2$ and $m^2+pe^2 =2a$.
\end{proof}
\end{prop}
$\bullet$ Next proposition tell us that there is a relationship between solutions $(b_{2})$ and Pythagorean triples such that difference of the two smaller sides is square.
\begin{prop}
If $pe^2=2a - m^2$, then $(m, e, u)$, is a solution $(b_{2})$ if and only if $(a-m^2, a, u)$, is a
Pythagorean triple.
\begin{proof}
$(a - m^2)^2 + a^2 = u^2$ if and only if $m^2(2a - m^2)= 2a^2 - u^2$
if and only if $[pe^2m^2=2a^2 - u^2$ and $pe^2=2a - m^2]$.
\end{proof}
\end{prop}

\begin{rem}
$(a, b, u), a<b$ is primitive Pythagorean triple if and only if there are coprime positive integers $s$ and $t$ such that
\begin{center}
$a=s^2 - t^2$ , $b=2st$\qquad or \qquad $a=2st$ , $b=s^2 - t^2$.
\end{center}
\end{rem}
\begin{lem}
If $a = s^2 - t^2 , b =2s t , (s , t)=1$ and $a - b= m^2$, then there are integers $x$ and $y$ such
that $s = 2x^2 + y^2 + 2xy , t = 2xy$ and $(2x , y)=1$.
\begin{proof}
We have $(s - t)^2 - 2t^2=m^2$. Let $k= s - t$, therefore $2t^2= k^2 - m^2$. From $(s , t)=1$ we have $(k , m)=1$.
So $m$ , $k$ are odd and $(k-m , k+ m)=2$. Consequently, there are coprime integers $x$ and $y$ such that
\begin{center}
 $k+m= 4x^2 , ~ k - m= 2y^2 \qquad\text {or}\qquad  k - m= 4x^2 , ~ k+m= 2y^2$.
\end{center}
Then $k= y^2+2x^2$ and $m= \pm(y^2 - 2x^2)$. Hence $t= \pm 2xy$ and $s=2x^2 + y^2 \pm 2xy$.
\end{proof}
\end{lem}
\begin{cor}
$a= 4x^4+y^4+ 4x^2 y^2+8x^3 y+4xy^3$
\begin{center}
and

$b= 8x^2 y^2+ 8x^3 y+4xy^3$.
\end{center}
\end{cor}
\begin{lem}
If $a = 2st, b = s^2 - t^2 , (s , t)=1$ and $a - b=m^2$ then there are integers $x$ and $y$ such that
$s = 2xy , t = 2x^2 + y^2+2xy$ and $(2x, y)=1$.
\begin{proof}
The proof is similar to the previous lemma by  letting $k= s + t$.
\end{proof}
\end{lem}

\begin{cor}
$b= - 4x^4-y^4-4x^2 y^2+8x^3 y+4xy^3$
\begin{center}
and

$a= - 8x^2 y^2+8x^3 y+4xy^3$.
\end{center}
\end{cor}

\begin{prop}
$(b_{2})$ has integer solution if and only if $p\square\in Imf_{1}\cup Imf_{2},$ where

$f_{1}(x, y) =~~~~~~ 4x^4 + y^4+12 x^2 y^2+16x^3 y+8xy^3, ~~~(2x, y)=1$,

$f_{2}(x, y) = -4x^4 -y^4-12x^2 y^2+16x^3 y+8xy^3, (2x, y)=1$.

\begin{proof}
$(b_{2})$ has integer solution if and only if $pe^2= a+b,~ a - b=m^2$ and $a^2+b^2= u^2$, then the above above corollaries yield the result.
\end{proof}
\end{prop}

\begin{exa}
$f_{1}(1,1)=41, f_{1}(-1,7)=137, f_{2}(1,1)=7$, hence for $p= 41, 137, 7$ equation $(b_{2})$ has an integer solution.
\end{exa}

\begin{rem}
As $y$ is odd we have $f_{1}(x,y)\equiv (2x^2+y)^2\equiv1\pmod 8$ and $f_{2}(x,y)\equiv -(2x^2+y)^2\equiv 7\pmod 8$.
\end{rem}

\begin{prop}
If $pe^2$ is $\beta_{-}$Pythagorean then $p$ is $\beta_{-}$Pythagorean.
\begin{proof}
If $pe^2$ is $\beta_{-}$Pythagorean if and only if there are integers $m, E$ and $u$ such that
$2u^2= m^4+(pe^2)^2 E^4 , (m , E) =1$.

Consequently $m$ and $e$ are odd. If $(m, e) = d,$ then $d$ is odd and also $d^2 \mid u$. There are integers $m_{0} ,u_{0 }$ and $e_{0}$
such that $e=e_{0}d,~ m=m_{0}d$ and $u=u_{0}d^2$. We know $(m_{0} , e_{0})=1$. Hence $m_{0}^4+p^2 (e_{0} E)^4=2u_{0}^2 , (m_{0} , e_{0} E)= 1$. Consequently $p$ is $\beta_{-}$Pythagorean.
\end{proof}
\end{prop}

\begin{cor}
$(b_{2})$ has an integer solution if and only if $p\in Imf_{1}\cup Imf_{2}$.

i.e., for solving the equation $(b_{2})$, we can choose $e=1$.
\end{cor}

\begin{exa} $17$ is not $\beta_{-}$Pythagorean because, there are not positive integers $a$ and $b$ such that $17=a+b,~ a-b=m^2~$ and $~a^2+b^2=u^2$.
\end{exa}

\begin{exa} $[41= 21+20,~ 21-20= \square~$ and $~21^2+20^2= \square]$, i.e., $f_{1}(1,1)=41$, therefore $41$ is $\beta_{-}$Pythagorean.
\end{exa}

$\bullet$ Now, we study the solvability of the homogeneous equation $(b_{3})$.
\begin{prop}
If $(b_{3})$ has an integer solution, then $p\square \in Imf_{3}$, where $f_{3}(x, y) =16x^4+y^4+24x^2y^2,  (2x, y)=1$.
\begin{proof}
If $(b_{3})$ has an integer solution, then

\begin{center}
$N^2= 2p(m^4 +e^4) ~, ~ gcd(m, 2pe)=gcd(e, 2p)=1$.
\end{center}

There is an integer number $u$ such that $2pu^2=m^4+e^4.$

So $(m^2-e^2)^2=2(pu^2-e^2m^2)$. Again there is an integer number $a$ such that $pu^2-e^2 m^2= 2a^2, m^2-e^2= 2a$, from which one gets integers $y$ and $t$ such that $m - e= 2y, m + e = 2t$ and $a= 2yt$. Therefore $m= t + y$ and $e= t - y$, where $y$ or $t$ is odd and the other is even.

We have $pu^2 = t^4 +y^4 +6y^2 t^2$. Suppose $t$ is even and $y$ is odd. There is an integers $x$ such that $t=2x$. This yields the result.
 \end{proof}
 \end{prop}
 \begin{exa}
$ f_{3}(1, 1)= 41,~ f_{3}(2, 1)= 353$, hence for $p= 41, ~353$, the equation $(b_{3})$ has an integer solution.
 \end{exa}

\begin{prop}
If $(b_{2})$ or $(b_{3})$ has an integer solution then equation $(a)$ has an integer solution.
\begin{proof}
If $(b_{2})$ has an integer solution, then $N^2=2(p^2e^4+m^4) , (me , 2)=1$. There is an integer number $u$ such that
$2u^2=p^2e^4+m^4$. Let $2c= pe^2+m^2, ~2d= pe^2-m^2$. We have $c+d= pe^2, ~c- d= m^2$. This implies that $c^2-d^2= pe^2m^2,~ c^2+d^2= u^2$. Hence $c^4 - d^4= p(emu)^2$.

If $(b_{3})$ has an integer solution, then $N^2= 2p(m^4+e^4) , (me , 2)=1$. Therefore, there is an integer number $u$ such that
$2pu^2=m^4+e^4$. Let $2c= e^2+m^2, ~2d= e^2 - m^2$. We have $c+d= e^2,~ c - d=m^2$. This implies that $c^2-d^2= e^2m^2, ~ c^2+d^2= pu^2$. Hence $c^4 - d^4= p(emu)^2$.
\end{proof}
\end{prop}

\begin{cor}
If $(b_{2})$ or $(b_{3})$ has an integer solution, then $p$ is congruent.
\begin{proof}
If $(b_{2})$ or $(b_{3})$ has an integer solution then
$\mid Im\overline\alpha\mid\geq 2$ and $\mid Im\alpha\mid =4$. As \eqref{eq2}, we have $r(E)\geq 1$. Consequently
$p$ is a congruent number.
\end{proof}
\end{cor}

\begin{mainthm}
$r(E)=2$ if and only if $p$ is $\alpha_{-}^{\pm}$Pythagorean and $\beta_{-}$Pythagorean.
\begin{proof}
$p$ is $\alpha_{-}^{\pm}$Pythagorean and $\beta_{-}$Pythagorean if and only if  $p\in Im \overline{\alpha}$ and $2 \in Im \overline{\alpha}$  if and only if  $(b_{1})$ and $(b_{2})$ have integer solutions if and only if $Im\overline{\alpha}= \lbrace 1, 2, p, 2p \rbrace$ and $Im{\alpha}=\lbrace \pm1, \pm p \rbrace$ if and only if $r(E)=2$.
\end{proof}
 \end {mainthm}

\begin{cor}
If $p\equiv1\pmod 8$ then $r(E)= 2$, whenever $p$ satisfies in the two following conditions:

(1) $\exists a , b \in \mathbb{N}$: $p= a+b~, ~ a - b= \square$ and $a^2+b^2= \square$.

(2) $\exists a , b \in \mathbb{N}$: $p= a^2+b^2~, ~ (a, 2b)=1$ and $(a \pm 2b)^2+b^2= \square$.
\end{cor}

\begin{exa}
 For $p=41$ we have $r(E)= 2$.
\end{exa}

$~$
\subsection*{Part (3)}
$~$

According to the 2-descent method, for the elliptic curves
\begin{equation}\label{eqC}
E: y^2=x^3 - 4p^2 x
\end{equation}
\ \ \ \ and
\begin{equation}\label{eqD}
\overline{E }: y^2=x^3 +16 p^2 x.
\end{equation}
We have respectively
${\{\pm1}\} \subseteq Im\alpha \subseteq {\{\pm 1, \pm 2, \pm p, \pm 2p }\}$ and
$ {\{1}\} \subseteq Im \overline{\alpha} \subseteq {\{1, 2, p, 2p}\}$.
Therefore, according to \eqref{eq2}, the maximum rank of \eqref{eqC} is $3$.
In fact, we showed that $Im \overline{\alpha} \subseteq {\{1, p}\}$. Hence the maximum rank is $2$.

The homogeneous equations of $E$ are
\begin{equation*}
\begin{array}{l}
(a_{1}):N^2=\pm(m^4-4pe^4),\quad  \gcd(m,2pe)=\gcd(e,p)=\gcd(N,em)=1,\\
\\
(a_{2}):N^2=\pm(2m^4-2p^2e^4),\quad \gcd(m,2pe)=\gcd(e,2)=\gcd(N,em)=1,\\
\\
(a_{3}):N^2=\pm(2pm^4-2pe^4),\quad  \gcd(m,2pe)=\gcd(e,2p)=\gcd(N,em)=1
\end{array}
\end{equation*}
and the homogeneous equations of $\overline{E}$ are
\begin{equation*}
\begin{array}{l}
(b_{1}): N^2= 2m^4 +8p^2e^4,\quad \gcd(m,2pe)=\gcd(e,2)=\gcd(N,em)=1,\\
\\
 (b_{2}): N^2= 2pm^4 +8pe^4,\quad \gcd(m,2pe)=\gcd(e,2p)=\gcd(N,em)=1,\\
 \\
 (b_{3}): N^2= pm^4 +16pe^4,\quad \gcd(m,2pe)=\gcd(e,p)=\gcd(N,em)=1.
\end{array}
 \end{equation*}
First we study the solvability of the above homogeneous equations.
\begin{prop}\label{p37}
 If $(a_{2})$ has an integer solution, then $p\equiv1, 3, 7 \pmod 8$.
 \end{prop}
\begin{proof}
If $(a_{2})$ has an integer solution, then $ N^2=\pm 2(m^4 - p^2e^4)$, $gcd (m,p)=1$.
There is an integer number $u$ such that $2u^2=\pm(m^4 - p^2 e^4)$.
So $2u^2 \equiv \pm m^4\pmod p$. Hence $(2um^{*^2})^2 \equiv \pm2\pmod p$. Implies $(\frac{\pm2}{p})=+1$, i.e., $p\equiv 1, 3, 7 \pmod 8$.
\end{proof}

\begin{prop}\label{p38}
$(b_{1})$ and  $(b_{2})$  do not  have any integer solutions.
\end{prop}
\begin{proof}
If $(b_{1})$ has an integer solutio,n then $ N^2= 2(m^4 +4p^2e^4)$, $\gcd (m, 2)=1$.
There is an integer number $u$ such that $2u^2= m^4+4p^2e^4$. So $m$ is an even number, which is a contradiction.

 If $(b_{2})$ has an integer solutio,n then $ N^2= 2p (m^4  + 4e^4)$, $\gcd (m, 2)=1$.

There is an integer number $u$ such that $2pu^2= m^4 + 4e^4$. So $m$ is an even number, which is a contradiction.
\end{proof}

\begin{prop}\label{p39}
If $(b_{3})$ has an integer solution, then $p\equiv1 \pmod 8$.
\end{prop}
\begin{proof}
If $(b_{3})$ has integer solution, then $ N^2= p(m^4 +16e^4)$, $\gcd (m, 2)=1$.
There is an integer number $u$ such that $pu^2= m^4+16e^4$.

We have $pu^2 \equiv 1 \pmod 8$. So $(u, 8)=1$. Therefore $u^2 \equiv 1 \pmod 8$. Consequently $p \equiv 1 \pmod 8$.
\end{proof}

\begin{cor}\label{c310}
The above propositions are summarized in the followings:
\begin{enumerate}
\item $\pm p \in Im\alpha$ if and only if $(a_{1})$ has an integer solution.

\item $\pm 2 \in Im\alpha$ if and only if $(a_{2})$ has an integer solution. In this case,

$p\equiv 1, 3, 7\pmod8$.

\item $\pm2p \in Im\alpha$ if and only if $(a_{3})$ has an integer solution.

\item $2\notin Im\overline{\alpha}$, i.e,. $(b_{1})$ does not have an integer solution.

\item $2p\notin Im\overline{\alpha}$,  i.e., $(b_{2})$ does not have an integer solution.

\item $p \in Im\overline{\alpha}$ if and only if $(b_{3})$ has an integer solution. In this case,

 $p \equiv 1\pmod 8$.
\end{enumerate}
\end{cor}
\begin{cor}\label{c311}
By using the previous corollary, we have:
\begin{enumerate}
\item Let $p\equiv 5 \pmod 8$.

$(b_{3})$ does not have integer solutions, then $Im\overline{\alpha}= \lbrace1\rbrace$.

 It is possible for $(a_{1})$ and $(a_{3})$ to have integer solutions. So $Im{\alpha}\subseteq \lbrace\pm 1, \pm p \rbrace$ or $ \lbrace\pm1, \pm 2p \rbrace $.

As \eqref{eq2}, implies $\mid Im\alpha\mid\geq4$.
Therefore $Im{\alpha}= \lbrace\pm 1, \pm p \rbrace$ or $\lbrace\pm1, \pm 2p \rbrace$.  Hence $r= 0$, i.e., 2p is not congruent.
\\
\item Let $p\equiv 3, 7\pmod 8$.

$(b_{3})$ does not have integer solutions, then $Im\overline{\alpha}= \lbrace1\rbrace$.

It is possible for $(a_{1})$, $(a_{2})$ and $(a_{3})$ to have integer solutions, then $Im{\alpha}\subseteq \lbrace\pm1, \pm 2, \pm p, \pm 2p \rbrace$.

$2p$ is congruent, (see\cite{mon}), so $r\geq1$.

The equation \eqref{eq2}, implies that $\mid Im\alpha\mid\geq8$. Therefore $Im{\alpha}\ = \lbrace\pm1, \pm 2, \pm p, \pm 2p \rbrace$. Hence $r = 1$.
\end{enumerate}
\end{cor}
$~$
\subsection*{part (4)}
$~$

According to the 2-descent method, for the elliptic curves
\begin{equation}\label{eqE}
E: y^2=x^3 - p^2q^2 x
\end{equation}
\ \ \ \ and
\begin{equation}\label{eqF}
\overline{E }: y^2=x^3 +4 p^2q^2 x.
\end{equation}
We have respectively
${\{\pm1}\} \subseteq Im\alpha \subseteq {\{\pm 1 ,\pm p,\pm q \pm pq}\}$ and
$ {\{1}\} \subseteq Im \overline{\alpha} \subseteq {\{1, 2, p, q, pq, 2p, 2q, 2pq}\}$.
Therefore, according to \eqref{eq2}, the maximum rank of \eqref{eqE} is 4. Moreover, the homogeneous equations of $E$ are
\begin{equation*}
\begin{array}{l}
(a_{1}): N^2=\pm( pm^4 - pq^2e^4),\quad \gcd(m, pqe)=\gcd(e, pm)=\gcd(N, me)=1,\\
\\
(a_{2}): N^2= \pm(qm^4 -p^2qe^4),\quad \gcd(m, pqe)=\gcd(e, qm)=\gcd(N, me)=1,\\
\\
(a_{3}): N^2= \pm(pqm^4 - pqe^4),\quad \gcd(m, pqe)=\gcd(e, pqm)=\gcd(N, me)=1,
\end{array}
\end{equation*}
and the homogeneous equations of $\overline{E}$ are
\begin{equation*}
\begin{array}{l}
(b_{1}): N^2= pm^4 +4pq^2e^4,\quad \gcd(m, 2pqe)=\gcd(e, pm)=\gcd(N, me)=1,\\
\\
(b_{2}): N^2= qm^4 +4p^2qe^4,\quad \gcd(m, 2pqe)=\gcd(e, qm)=\gcd(N, me)=1,\\
\\
(b_{3}): N^2= 2m^4 +2p^2q^2e^4,\quad \gcd(m, 2pqe)=\gcd(e, 2m)=\gcd(N, me)=1,\\
\\
(b_{4}): N^2= 2pm^4 +2pq^2e^4,\quad \gcd(m, 2pqe)=\gcd(e, 2pm)=\gcd(N, me)=1,\\
\\
(b_{5}): N^2= 2qm^4 +2p^2qe^4,\quad \gcd(m, 2pqe)=\gcd(e, 2qm)=\gcd(N, me)=1,\\
\\
(b_{6}): N^2= pqm^4 +4pqe^4,\quad \gcd(m, 2pqe)=\gcd(e, pqm)=\gcd(N, me)=1,\\
\\
(b_{7}): N^2= 2pqm^4 +2pqe^4,\quad \gcd(m, 2pqe)=\gcd(e, 2pqm)=\gcd(N, me)=1.
\end{array}
\end{equation*}
Equations $(a_{1}), (a_{2})$ and $(b_{1}), (b_{2})$ and $(b_{4}), (b_{5})$ are the same.

 First we study the solvability of these homogeneous equations.
\begin{prop}\label{312}
If $(a_1)$ has an integer solution, then $[p$ or $q \equiv 1, 3, 7 \pmod 8$ and $(\frac{-q}{p})=+1]$ or $[p$ or $q \equiv 1, 7 \pmod 8$ and $(\frac{q}{p})=+1]$.
\begin{proof}
If $(a_1)$ has an integer solution, then $ N^2 =\pm p(m^4 -q^2e^4)$, $\gcd(m,Nepq)=1$. There is an integer number $u$ such that $\pm pu^2= m^4- q^2e^4$.

If $q\mid u$, then $q\mid m$. From $N=pu$, we have $\gcd(N,m)\neq 1$, a contradiction. So $\gcd(q,u)=1$.

There are positive integers $u_{1}$, $u_{2}$ and  a square free $t$  such that
\begin{equation*}
m^2+ qe^2=pu_{1}^2t\qquad \text{and}\qquad m^2- qe^2=\pm u_{2}^2t,
\end{equation*}
or
\begin{equation*}
m^2- qe^2=\pm pu_{1}^2t\qquad \text{and}\qquad m^2+qe^2=u_{2}^2t.
\end{equation*}

We have $2m^2= \pm t(pu_{1}^2 \pm u_{2}^2)$ and $2qe^2= \mp t(pu_{1}^2 \mp u_{2}^2)$. So $t \mid 2m^2$ and
$t\mid 2qe^2$.

Suppose $t$ is odd number. As  $t\mid m^2$, we have $\gcd(m, t)\neq 1$. So $\gcd(t, e)=1$. From $t\mid qe^2$ we have $t\mid q$. Because $u=u_{1}u_{2}t$ and $gcd(q, u)=1$, we have $t=1$, which is a contradiction. Therefore $t=1$.
\begin{enumerate}
\item If $m^2+ qe^2= pu_{1}^2$ and $m^2- qe^2=\pm u_{2}^2$ then $2m^2= pu_{1}^2 \pm u _{2}^2$.

From $m^2+ qe^2= pu_{1}^2$ we have, $(\frac{-q}{p})=+1$.

From $2m^2= pu_{1}^2 \pm u _{2}^2$ we have, $(\frac{\pm2}{p})=+1$, i.e., $p\equiv 1, 3, 7\pmod 8$.

\item If $m^2- qe^2=\pm pu_{1}^2$ and $m^2+qe^2= u_{2}^2$ then $2m^2= \pm pu_{1}^2+u _{2}^2$.

From $m^2- qe^2=\pm pu_{1}^2$, we have $(\frac{q}{p})=+1$.

From $2m^2= \pm pu_{1}^2+u _{2}^2$, we have $(\frac{2}{p})=+1$, i.e., $p\equiv 1, 7\pmod 8$.

Now, we suppose $t$ is even. If $t=2k$, where $k\neq 1$ is odd number, then $m^2= \pm k(pu_{1}^2 \pm u_{2}^2)$ and $qe^2= \mp k(pu_{1}^2 \mp u_{2}^2)$. So $k \mid m^2$ and $k\mid qe^2$, as above we have a contradiction. Hence $t=2$.

\item If $m^2 + qe^2= 2pu_{1}^2$ and $m^2 - qe^2=\pm 2u_{2}^2$ then $qe^2 = p u_{1}^2 \pm u_{2}^2$.

From $m^2 +qe^2=2pu_{1}^2$, we have $(\frac{-q}{p})=+1$.

$m^2-qe^2=\pm 2u_{2}^2$, so $(\frac{\pm2}{q})=+1$, i.e., $q\equiv1, 3, 7\pmod 8$.

\item If $m^2 - qe^2= \pm 2pu_{1}^2$ and $m^2 + qe^2=2u_{2}^2$ then $qe^2 = u_{2}^2 \mp p u_{1}^2$.
\end{enumerate}

From $m^2+qe^2=2u_{2}^2$, we have $(\frac{2}{q})=+1$, i.e., $q\equiv 1, 7\pmod 8$.

From $qm^2 = u_{2}^2 \pm p u_{1}^2$, we have $(\frac{q}{p})=+1$.
\end{proof}
\end{prop}

\begin{prop}
If $(b_{1})$ has an integer solution, then $p\equiv1 \pmod 4$ and $(\frac{p}{q})=+1$.
\begin{proof}
If $(b_{1})$ has an integer solution, then $ N^2= p(m^4 +4q^2e^4)$, $gcd(m, pq)=1$. There is an integer number $u$ such that $pu^2= m^4 +4q^2e^4$.

$m^4 \equiv-4q^2e^4 \pmod p$, so $(2qe^2 m^{\ast^2})^2\equiv-1 \pmod p$. Hence $(\frac{-1}{p})=+1$, i.e., $p\equiv1 \pmod 4$.

$pu^2\equiv m^4 \pmod q$, so $(m^{*^2}p u)^2 \equiv p\pmod q$, i.e., $(\frac{p}{q})=+1$.
\end{proof}
\end{prop}

\begin{prop}\label{p314}
If $(b_{3})$ has an integer solution, then $p, q\equiv\pm 1 \pmod 8$.
\end{prop}
\begin{proof}
If $(b_{3})$ has an integer solutio,n then $N^2= 2(m^4 +p^2q^2e^4)$, $\gcd(m, pq)=1$. There is an integer number $u$ such that $2 u^2= m^4 +p^2 q^2e^4$.

$2 u^2\equiv m^4\pmod p$, so $(m^{*^2}2u)^2\equiv2\pmod p$. This implies that $(\frac{2}{p})=+1$, i.e., $p\equiv\pm1 \pmod8$.

$2 u^2\equiv m^4\pmod q$, so $(m^{*^2}2u)^2\equiv2\pmod q$. This implies that $(\frac{2}{q})=+1$, i.e., $q\equiv\pm1 \pmod8$.
\end{proof}

\begin{prop}\label{p315}
If $(b_{4})$ has an integer solution, then $p\equiv1\pmod 4$ and $(\frac{2p}{q})=+1$.
\end{prop}
\begin{proof}
If $(b_{4})$ has an integer solution, then $ N^2= 2p(m^4 +q^2e^4)$, $\gcd(m, pq)=1$.
 There is an integer $u$ such that $2pu^2=m^4+q^2e^4$.

$2p u^2\equiv m^4 \pmod q$, so $(m^{*^2}2p u)^2\equiv 2p \pmod q$, i.e., $(\frac{2p}{q})=+1$.

$-m^4\equiv q^2e^4\pmod p$, so $(qe^2{m^\ast}^2)^2\equiv -1 \pmod p$. Hence $(\frac{-1}{p})=+1$, i.e., $p\equiv1\pmod4$.
\end{proof}

\begin{prop}\label{p316}
If $(b_{6})$ has an integer solution. then $ p, q\equiv1 \pmod 4$.
\end{prop}
\begin{proof}
If $(b_{6})$ has an integer solution, then $ N^2= pq (m^4 + 4e^4)$, $\gcd(m, 2pq)=1$. There is an integer number $u$ such that $pq u^2= m^4 +4e^4$.

$-m^4\equiv 4e^4 \pmod p$, so $(2e^2{m^\ast}^2)^2\equiv -1 \pmod p$. Hence $(\frac{-1}{p})=+1$, i.e., $p\equiv1\pmod4$.

$-m^4\equiv 4e^4 \pmod q$, so $(2e^2{m^\ast}^2)^2\equiv -1 \pmod q$. Hence $(\frac{-1}{q})=+1$, i.e., $q\equiv1\pmod4$.
\end{proof}

Finally we verify $(b_{7})$. First we prove the following lemma.
\begin{rem}\label{r317}
If $b$ is a quadratic residue, then obviously $b^*$ is a quadratic residue. Moreover if $-1$ is a quadratic residue then $-b$ is a quadratic residue too.
\end{rem}

\begin{rem}\label{r318}
$-1$ is a quadratic residue if and only if $ p \equiv 1\pmod4$, i.e., $p \equiv 1, 5 \pmod 8$.
\end{rem}

\begin{lem}\label{l319}
If there is an integer number $x$ such that $x^4 \equiv -1 \pmod p$ then $p \equiv 1 \pmod 8$.
\end{lem}
\begin{proof}
$\frac{p-1}{2}$ number are quadratic residues and the same number of non-residues.

If $p=8k+5$, then there are $4k+2$ quadratic residues.

\begin{center}
$x^2\equiv\pm b \pmod 8 \Longleftrightarrow x^4\equiv b^2 \pmod 8$.
\end{center}

Therefore, there are $2k+1$ residues of degree $4$. So there are $2k$ residues of degree $4$ except 1.
If $b\not=\pm1$ is a residue, then $b^\ast\not=b$ and $b^\ast$ is a residue as well. Moreover only $\pm 1$ equal to their inverse in every mod. Putting these together, one can get that $-1$ is not residue of degree $4$.
\end{proof}

\begin{prop}\label{p320}
If $(b_{7})$ has an integer solution, then $p, q\equiv1\pmod 8$.
\end{prop}
\begin{proof}
If $(b_{7})$ has an integer solution, then $N^2=2pq(m^4 +e^4)$, $\gcd(m, pq)=1$. There is integer number $u$ such that $2pq u^2= m^4 +e^4$.

$-m^4\equiv e^4 \pmod p$, so $(e{m^\ast})^4\equiv -1 \pmod p$, i.e., $p\equiv1\pmod8$.

$-m^4\equiv e^4 \pmod q$, so $(e{m^\ast})^4\equiv -1 \pmod q$, i.e., $q\equiv1\pmod8$.
\end{proof}

\begin{cor}\label{c321}
The above propositions are summarized in the following statements.
\begin{enumerate}
\item $p\in Im\alpha$ if and only if $(a_{1})$ has an integer solution. In this case,

$p$ or $q \equiv 1, 3, 7 \pmod 8$ and $(\frac{-q}{p})=+1$

 or

 $p$ or $q \equiv 1, 7 \pmod 8$ and $(\frac{q}{p})=+1$.

\item $q\in Im\alpha$ if and only if $(a_{2})$ has an integer solution. In this case,

$p$ or $q \equiv 1, 3, 7 \pmod 8$ and $(\frac{-p}{q})=+1$

or

$p$ or $q \equiv 1, 7 \pmod 8$ and $(\frac{p}{q})=+1$.

\item $pq \in Im\alpha$ if and only if $ (a_{3})$ has an integer solution.

\item $p\in Im\overline{\alpha}$ if and only if $(b_{1})$ has an integer solution. In this case,

$p\equiv1, 5\pmod 8$ and $(\frac{p}{q})=+1$.

\item $q\in Im\overline{\alpha}$ if and only if $ (b_{2})$ has an integer solution. In this case,

$q\equiv1, 5\pmod 8$ and $(\frac{q}{p})=+1$.

\item $2\in Im\overline{\alpha}$ if and only if $(b_{3})$ has an integer solution. In this case,

$p, q\equiv 1, 7\pmod 8$.

\item $2p\in Im\overline{\alpha}$ if and only if $(b_{4})$ has an integer solution. In this case,

$p\equiv1, 5, ~ q\equiv 1, 7\pmod 8$ and $(\frac{p}{q})=+1$

or

$p\equiv1, 5, ~q\equiv 3, 5\pmod 8$ and $(\frac{p}{q})=-1$.

\item $2q\in Im\overline{\alpha}$ if and only if $(b_{5})$ has an integer solution' In this case,

$q\equiv1, 5,~ p\equiv 1, 7 \pmod 8$ and $(\frac{q}{p})=+1$

 or

$q\equiv1, 5,~p\equiv 3, 5 \pmod 8$ and $(\frac{q}{p})=-1$.

\item $pq\in Im\overline{\alpha}$ if and only if $(b_{6})$ has an integer solution. In this case,

$p, q\equiv1, 5\pmod8$.

\item $2pq\in Im\overline{\alpha}$ if and only if $(b_{7})$ has an integer solution. In this case,

$p, q\equiv 1\pmod8$.
\end{enumerate}
\end{cor}

\begin{cor}\label{c322}
By the previous corollary, we have:
\begin{enumerate}
\item Let $ p, q\equiv3\pmod8$.

$(b_{i}), 1\leq i \leq 7$, does not have an integer solution. So $Im\overline{\alpha}=\lbrace1\rbrace$.

As \eqref{eq2}, then $\mid Im\alpha\mid\geq4$.

It is possible for $(a_{1}), (a_{2})$ and $(a_{3})$ to have integer solutions. If $(a_{1})$ and $(a_{2})$  have integer solutions, then $(\frac{q}{p})=-1$ and $(\frac{p}{q})=-1$. By using quadratic reciprocity, which is a contradiction. Therefore $Im{\alpha}\subseteq \lbrace \pm 1, \pm p \rbrace,  \lbrace \pm 1, \pm q \rbrace$ or $\lbrace \pm 1, \pm pq \rbrace$.

Consequently, $Im{\alpha}= \lbrace \pm 1, \pm p \rbrace , \lbrace \pm 1, \pm q \rbrace$ or $\lbrace \pm 1, \pm pq \rbrace$.

Hence $r= 0$, i.e., $pq$ is not congruent.
\\
\item Let $ p\equiv3,~ q\equiv5 \pmod8$.

It is possible for $(b_{2})$ or $(b_{5})$ to have an integer solution. Therefore $Im\overline{\alpha}\subseteq\lbrace1, q \rbrace$ or $\lbrace 1, 2q \rbrace$.

It is possible for $(a_{3})$ to have an integer solution, then $Im{\alpha} \subset \lbrace \pm 1, \pm pq \rbrace$.
As $pq$ is congruent number, (see\cite{mon}). So $ r\geq1$.

As \eqref{eq2}, implies that $Im{\alpha} = \lbrace \pm 1, \pm pq \rbrace$ and $Im\overline{\alpha}=\lbrace1, q \rbrace$ or $\lbrace 1, 2q \rbrace$. Hence $r= 1$.
\\
\item Let $ p\equiv 3,~ q\equiv7\pmod8$.

 $(b_{i}), 1\leq i \leq 7$, does not have an integer solution. So $Im\overline{\alpha}=\lbrace1\rbrace$.

It is possible for $(a_{1}), (a_{2})$ and $(a_{3})$ to have integer solutions. So $Im{\alpha} \subseteq \lbrace \pm 1, \pm p, \pm q, \pm pq \rbrace$.

As $pq$ is congruent number, (see\cite{mon}), so $ r\geq1$.

As \eqref{eq2},implies that  $Im{\alpha} = \lbrace \pm 1, \pm p, \pm q, \pm pq \rbrace$. Hence $r= 1$.
\\
\item Let $ p\equiv1, ~q\equiv5 \pmod8$ and $(\frac{p}{q})=-1$.

It is possible for $(a_{3})$ to have an integer solution. Therefore, $Im{\alpha} \subseteq \lbrace \pm 1, \pm pq \rbrace$.

It is possible for $(b_{4})$ or $(b_{6})$ to have an integer solution. Therefore, $Im\overline{\alpha}\subseteq\lbrace1, 2p \rbrace$ or $\lbrace 1, pq \rbrace$.

As $pq$ is congruent number, (see\cite{mon}), we have $r\geq1$.

As \eqref{eq2}, Implies that $Im\overline{\alpha}=\lbrace1, 2p \rbrace$ or $\lbrace 1, pq \rbrace$ and
$Im{\alpha}=\lbrace\pm1, \pm pq \rbrace$. Hence $r=1$.
\\
\item Let $p\equiv1,~ q\equiv7 \pmod8$ and $(\frac{p}{q})=-1$.

It is possible for $(a_{2})$ or $(a_{3})$ to have an integer solution. Therefore, $Im{\alpha} \subseteq \lbrace \pm1, \pm q \rbrace$ or $\lbrace \pm1, \pm pq \rbrace$.

It is possible for $(b_{3})$ to have an integer solution. Then $Im\overline{\alpha}\subseteq\lbrace1, 2 \rbrace$.

As $pq$ is a congruent number, (see\cite{mon}), we have $r\geq1$.

As \eqref{eq2}, Implies that  $Im{\alpha}= \lbrace \pm1, \pm q \rbrace$ or $\lbrace \pm1, \pm pq \rbrace$ and $Im\overline{\alpha}=\lbrace1, 2 \rbrace$. Hence $r=1$.
\end{enumerate}
\end{cor}
$~$
\subsection*{part (5)}
$~$

According to the 2-descent method, for the elliptic curves
\begin{equation}\label{eq3}
E: y^2=x^3 - 4p^2q^2 x
\end{equation}
\ \ \ \ and
\begin{equation}\label{eq30}
\overline{E }: y^2=x^3 +16 p^2q^2 x.
\end{equation}
We have respectively
${\{\pm1}\} \subseteq Im\alpha \subseteq {\{\pm 1,\ \pm2,\ \pm p,\ \pm q,\ \pm pq,\ \pm2p,\ \pm2q,\ \pm2pq}\}$ and
%\vspace*{3pt}
%$$ then
$ {\{1}\} \subseteq Im \overline{\alpha} \subseteq {\{1, 2, p, q, pq, 2p, 2q, 2pq}\}$.
Therefore, according to \eqref{eq2},the maximom rank of \eqref{eq3} is 5. Moreover, the homogeneous equations of $E$ are

\begin{equation*}
\begin{array}{l}
(a_{1}): N^2=\pm(2m^4 -2p^2q^2e^4),\quad \gcd(m, 2pqe)= \gcd(e, 2m)= \gcd(N, me)=1,\\
\\
(a_{2}): N^2=\pm(pm^4 -4pq^2e^4),\quad \gcd(m, 2pqe)= \gcd(e, pm)= \gcd(N, me)=1 ,\\
\\
(a_{3}): N^2= \pm(qm^4 -4p^2qe^4),\quad  \gcd(m, 2pqe)= \gcd(e, qm)= \gcd(N, me)=1,\\
\\
(a_{4}): N^2= \pm(2pm^4 -2pq^2e^4),\quad \gcd(m, 2pqe)= \gcd(e, 2pm)= \gcd(N, me)=1,\\
\\
(a_{5}): N^2=\pm(2qm^4 -2p^2qe^4),\quad \gcd(m, 2pqe)= \gcd(e, 2qm)= \gcd(N, me)=1,\\
\\
(a_{6}): N^2=\pm(pqm^4 -4pqe^4),\quad \gcd(m, 2pqe)= \gcd(e, pqm)= \gcd(N, me)=1,\\
\\
(a_{7}): N^2=\pm(2 pqm^4 -2pqe^4),\quad \gcd(m, 2pqe)= \gcd(e, 2pqm)= \gcd(N, me)=1,
\end{array}
\end{equation*}
and the homogeneous equations of $\overline{E}$ are
\begin{equation*}
\begin{array}{l}
(b_{1}): N^2= 2m^4 +8p^2q^2e^4,\quad  \gcd(m, 2pqe)= \gcd(e, 2m)= \gcd(N, me)=1,\\
\\
(b_{2}): N^2= pm^4 +16pq^2e^4,\quad \gcd(m, 2pqe)= \gcd(e, pm)= \gcd(N, me)=1,\\
\\
(b_{3}): N^2= qm^4 +16p^2qe^4,\quad \gcd(m, 2pqe)= \gcd(e, qm)= \gcd(N, me)=1,\\
\\
(b_{4}): N^2= 2pm^4 +8pq^2e^4,\quad \gcd(m, 2pqe)= \gcd(e, 2pm)= \gcd(N, me)=1,\\
\\
(b_{5}): N^2= 2qm^4 +8p^2qe^4,\quad \gcd(m, 2pqe)= \gcd(e, 2qm)= \gcd(N, me)=1,\\
\\
(b_{6}): N^2= pqm^4 +16pqe^4,\quad \gcd(m, 2pqe)= \gcd(e, pqm)= \gcd(N, me)=1,\\
\\
(b_{7}): N^2= 2pqm^4 +8pqe^4,\quad \gcd(m, 2pqe)= \gcd(e, 2pqm)= \gcd(N, me)=1.
\end{array}
\end{equation*}
\\
The equations $(a_{2}), (a_{3})$ and $(a_{4}), (a_{5})$ and $(b_{2}), (b_{3})$ and also $(b_{4}), (b_{5})$ are similar.

First we study the solvability of the above homogeneous equations.
\begin{prop}\label{p359}
If $(a_{1})$ has an integer solution, then $p, q\equiv 1, 3, 7\pmod 8$.
\end{prop}
\begin{proof}
If $(a_{1})$ has an integer solution, then $ N^2=\pm 2(m^4 -p^2q^2e^4)$, $gcd(m, Npq)=gcd(me, 2)=1$. There is an integer number $u$ such that $\pm2u^2= m^4 - p^2q^2 e^4$. There are  integer numbers $u _{1}$, $u _{2}$ and an odd number $t$ such that
\begin{equation*}
m^2-pqe^2=\pm2u_{2}^2t \qquad \text{and}\qquad m^2+pqe^2=4u_{1}^2t
\end{equation*}
or
\begin{equation*}
m^2-pqe^2=\pm4u_{1}^2t \qquad \text{and}\qquad m^2+pqe^2=2u_{2}^2t.
\end{equation*}

We have $m^2=\pm t(2u_{1}^2 \pm u_{2}^2)$. So $t\mid m^2$. We know $N=2u= 2u_{1}u_{2}t$ and $(N, m)=1$. Hence $t=1$.

If $p\mid u_{1}$, from $m^2\pm pqe^2=\pm 4u_{1}^2$, we have $p\mid m$, which is a contradiction. Consequently $(p, u_{1})=1$.

We have $pqe^2= \pm(2u_{1}^2 \pm u_{2}^2)$. Implies $\pm 2u_{1}^2 \equiv u_{2}^2 \pmod p$.
Therefore $\pm 2 \equiv (u_{2} u_{1}^\ast)^2\pmod p$. Consequently $(\frac{\pm 2}{p})=+1$, i.e., $p\equiv1, 3, 7 \pmod 8$.

The Proof for $q$ is similar.
\end{proof}

\begin{prop}\label{p360}
If $(a_{2})$ has an integer solutio,n then $(\frac{\pm p}{q})=+1 =(\frac{\pm 2q}{p})$.
\end{prop}
\begin{proof}

If $(a_{2})$ has an integer solution, then $ N^2=\pm p(m^4 -4q^2e^4)$, $\gcd(m, Npq)=\gcd(e,p)=1$. There is an integer number $u$ such that $\pm pu^2= m^4 -4q^2 e^4$. So, there are positive integer numbers $u _{1}$, $u _{2}$ and $t$ such that,
\begin{center}
$m^2-2qe^2=\pm pu_{1}^2t \qquad \text{and}\qquad m^2+2qe^2=u_{2}^2t$
\end{center}
or
\begin{center}
$m^2-2qe^2=\pm u_{2}^2t \qquad \text{and}\qquad m^2+2qe^2=pu_{1}^2t$.
\end{center}

We have $2m^2= \pm t(pu_{1}^2 \pm u_{2}^2)$. So $t \mid 2m^2$.
\begin{enumerate}
\item suppose $t$ is odd.

We have $t \mid m^2$. Because $N=pu=pu_{1}u_{2}t$ and $\gcd(N, m)=1$, consequently $t=1$.

From $m^2\equiv\pm2qe^2 \pmod p$, we have $(me^\ast)^2\equiv\pm2q \pmod p$. So $(\frac{\pm 2q}{p})=+1.$

From $m^2\equiv \pm p u_{1}^2 \pmod q$, we have $(pm^*u_{1})^2\equiv \pm p\pmod q$. Consequently $(\frac{\pm p}{q})=+1$.

\item Now, suppose $t$ is even.

If $t=2k$, where $k$ is an odd integer then $m^2= \pm k(pu_{1}^2 \pm u_{2}^2)$. So $k \mid m^2$, similarly we have $k=1$. Hence $t=2$.
\end{enumerate}
From $m^2\pm2qe^2=\pm 2pu_{1}^2$, we get $m^2\equiv \pm 2pu_{1}^2 \pmod q$. Hence $(m^*2pu_{1})^2\equiv \pm 2p \pmod q$.
Consequently $(\frac{\pm2p}{q})=+1$.

From $m^2 \pm 2qe^2=\pm 2u_{2}^2$, we get $m^2\equiv \pm 2u_{2}^2 \pmod q$. Hence $\pm 2\equiv (2m^*u_{2})^2 \pmod q$. Consequently $(\frac{\pm 2}{q})=+1$.

From $m^2\pm 2qe^2=\pm 2pu_{1}^2$, we get $m^2\equiv \pm 2qe^2 \pmod p$. Hence $(me^*)^2 \equiv \pm 2q\pmod q$. Consequently $(\frac{\pm 2q}{p})=+1$.
\end{proof}

\begin{cor}\label{c326}
If $p, q \equiv 5\pmod 8$, then $(a_{2})$ does not have an integer solution.
\end{cor}
\begin{proof}

If $(a_{2})$ has an integer solution, then $(\frac{ \pm p}{q})=+1$ and $(\frac{\pm 2q}{p})=+1$.

$(\frac{2q}{p})=+1$ if and only if $[p\equiv\ 1, 7 \pmod 8$ and $(\frac{q}{p})=+1]$ or $[p\equiv\ 3, 5 \pmod 8$ and $(\frac{q}{p})=-1]$.

$(\frac{-2q}{p})=+1$ if and only if $[p\equiv\ 1, 3 \pmod 8$ and $(\frac{q}{p})=+1]$ or $[p\equiv\ 5, 7 \pmod 8$ and $(\frac{q}{p})=-1]$.

If for $p, q\equiv\ 5 \pmod 8$ equation $(a_{2})$ has an integer solution, then $(\frac{q}{p})=-1$ and $(\frac{\pm p}{q})= (\frac{p}{q})=+1$. By the quadratic reciprocity, this is a contradiction.
\end{proof}

\begin{prop}\label{p327}
If $(a_{4})$ has an integer solutio,n then $(\frac{q}{p})=+1=(\frac{\pm 2p}{q})$.
\end{prop}
\begin{proof}
If $(a_{4})$ has an integer solution, then $ N^2= \pm 2p(m^4 -q^2e^4)$, $\gcd(m, 2Npqe)=1$.
There is an integer number $u$ such that $\pm 2pu^2= m^4 -q^2 e^4$. So, there are integer numbers $u _{1}$, $u _{2}$ and a free square odd number $t$ such that
\begin{center}
$m^2 + qe^2=2pu_{1}^2t \qquad \text{and}\qquad m^2 - qe^2=\pm 4u_{2}^2t$.
\end{center}
or
\begin{center}
$m^2 - qe^2= \pm 4pu_{1}^2t \qquad \text{and}\qquad m^2 + qe^2=2u_{2}^2t$.
\end{center}

We have $m^2=t(pu_{1}^2\pm 2u_{2}^2)$ or $m^2=t(\pm 2pu_{1}^2+u_{2}^2)$. So $t\mid m^2$. From $N=2pu=4pu_{1}u_{2}t$ and $(N, m)=1$, implies $t=1$.

\begin {enumerate}
\item In case one, we have $m^2 \equiv2pu_{1}^2\pmod q$. So $2p \equiv (2pu_{1}m^*)^2\pmod q$. Consequently $(\frac{2p}{q})=+1$.

From $m^2 \equiv \pm qe^2\pmod p$, we have $\pm q \equiv (qe^2m^*)^2 \pmod p$. consequently $(\frac{q}{p})=+1$.

\item In case two, we have $m^2 \equiv \pm 4pu_{1}^2 \pmod q$, implies $\pm p \equiv (2m^*pu_{1} )^2 \pmod q$. Consequently $(\frac{\pm p}{q})=+1$.
\end {enumerate}

$m^2 \equiv 2u_{2}^2 \pmod q$, so $2 \equiv (2m^*u_{2} )^2 \pmod q$, i.e., $(\frac{2}{q})=+1$.

$m^2 \equiv qe^2 \pmod p$, so $(me^*)^2 \equiv q \pmod p$, i.e., $(\frac{q}{p})=+1$.
\end{proof}
\begin{cor}\label{c328}
If $p, q \equiv 5\pmod 8$, then $(a_{4})$ does not have any  integer solution.
\end{cor}
\begin{proof}
If $(a_{4})$ has an integer solution, then $ (\frac{p}{q})=+1$, and $(\frac{\pm 2q}{p})=+1$.

$(\frac{2q}{p})=+1$ if and only if $[p\equiv\ 1, 7\pmod 8$ and $(\frac{q}{p})=+1]$ or $[p\equiv\ 3, 5\pmod 8$ and $(\frac{q}{p})=-1]$.
$(\frac{-2q}{p})=+1$ if and only if $[p\equiv\ 1, 3\pmod 8$ and $(\frac{q}{p})=+1]$ or $[p\equiv\ 5, 7\pmod 8$ and $(\frac{q}{p})=-1]$. If for $p, q \equiv 5\pmod 8$, equation $(a_{4})$ has integer solution then $(\frac{q}{p})=-1$ and $(\frac{p}{q})=+1$. By using quadratic reciprocity, which is a contradiction.
\end{proof}
\begin{prop}\label{p364}
If $(a_{6})$ has integer solution then $p, q\equiv 1, 3, 7 \pmod 8$.
\end{prop}
\begin{proof}

If $(a_{6})$ has an integer solution, then $ N^2=\pm pq(m^4 - 4e^4)$, $\gcd(e, pq)=1$. There is an integer number $u$ such that $\pm pqu^2= m^4 - 4 e^4$. So, there are positive integers $u _{1}$, $u _{2}$ and  $t$ such that
\begin{equation*}
\begin{array}{l}
m^2+ 2e^2=pu_{1}^2 t \qquad \text{and}\qquad m^2- 2e^2=\pm qu_{2}^2 t,\\
\\
m^2+ 2e^2=qu_{1}^2 t \qquad \text{and}\qquad m^2- 2e^2=\pm pu_{2}^2 t,\\
\\
m^2+ 2e^2=pqu_{1}^2 t \qquad \text{and}\qquad m^2- 2e^2=\pm u_{2}^2 t,\\
\\
m^2- 2e^2=\pm pqu_{1}^2 t \qquad \text{and}\qquad m^2+2e^2= u_{2}^2 t.
\end{array}
\end{equation*}

We have $m^2\equiv\pm 2e^2 \pmod p$. Hence $(me^*)^2\equiv\pm 2\pmod p$. So $(\frac{\pm 2}{p})=+1$, consequently $p\equiv 1, 3, 7 ~(mod~8)$.

We have $m^2\equiv\pm 2e^2 \pmod q$. Hence $(me^*)^2\equiv\pm 2 \pmod q$. So $(\frac{\pm 2}{q})=+1$, consequently $q\equiv 1, 3, 7 ~(mod~8)$.
\end{proof}

\begin{prop}
$(b_{1}), (b_{4}), (b_{5})$ and $(b_{7})$ does not have any integer solution.
\end{prop}
\begin{proof}
If $(b_{1})$ has an integer solutio,n then $N^2=2(m^4+4p^2 q^2 e^4), (m, 2)=1$. There is an integer number $u$ such that $ 2u^2=m^4+4p^2 q^2 e^4 $. Hence $m$ is even, which is a contradiction.

The proof for $(b_{4 }), (b_{5} )$ and $(b_{7})$ are similar.
\end{proof}

\begin{prop}\label{p330}
If $(b_{2})$ has an integer solution, then $p\equiv 1 \pmod 8$ and $(\frac{p}{q})=+1$.
\end{prop}
\begin{proof}
If $(b_{2})$ has an integer solutio,n then $N^2=p(m^4+16q^2 e^4 )$, $\gcd(m, 2pq)=1$. There is an integer $u$ such that $pu^2=m^4+16q^2 e^4$.

Hence $pu^2\equiv m^4 \pmod q$. So $(pum^{*^2})^2\equiv p \pmod q$. Consequently $(\frac{p}{q})=+1$.

We have $m^4\equiv -16q^2 e^4 \pmod p$. So $(4qe^2m^{*^2})^2\equiv -1 \pmod p$, so $(\frac{-1}{p})=+1$, i.e., $p\equiv1 \pmod 4$. So $p\equiv1, 5 \pmod 8$.

 If $p\equiv 5 \pmod 8$, then $5u^2\equiv 1\pmod8$. Hence $u^2\equiv 5\pmod8$, which is a contradiction. Consequently $p\equiv1 \pmod 8$.
\end{proof}
\begin{prop}\label{p331}
If $(b_{6})$ has an integer solution, then $p, q \equiv 1 \pmod 8$.
\end{prop}
\begin{proof}
If $(b_{6})$ has an integer solution, then $N^2=pq(m^4+16e^4 )$, $\gcd(m, pq)=1$. There is an integer $u$ that $pqu^2=m^4+16e^4$.

From $m^4\equiv -16e^4 \pmod p$, we have $(2e{m^*})^4\equiv -1 \pmod p$. Consequently $p\equiv 1\pmod 8$.

From $m^4\equiv -16e^4 \pmod q$, we have $(2e{m^*})^4\equiv -1 \pmod q$. Consequently $q\equiv 1\pmod 8$.
\end{proof}

\begin{cor}
The above results are summarized in the followings:
\begin{enumerate}
\item $2 \in Im \alpha$ if and only if $(a_{1})$ has an integer solution. In this case,

$p, q \equiv 1, 3, 7\mod8$.

\item $p \in Im\alpha$ if and only if $(a_{2})$ has an integer solution. In this case,

$(\frac{\pm p}{q})=+1= (\frac{\pm 2q}{p})$.

\item $q \in Im\alpha$ if and only if $(a_{3})$ has an integer solution. In this case,

$(\frac{\pm q}{p})=+1= (\frac{\pm 2p}{q})$.

\item $2p \in Im\alpha$ if and only if $(a_{4})$ has an integer solution. In this case,

$(\frac{q}{p})=+1=(\frac{\pm 2p}{q})$.

\item $2q \in Im\alpha$ if and only if $(a_{5})$ has an integer solution. In this case,

$(\frac{ p}{q})=+1=(\frac{\pm 2q}{p})$.

\item $pq \in Im\alpha$ if and only if $(a_{6})$ has an integer solution. In this case,

 $p, q\equiv 1, 3, 7\pmod8$.

\item $2pq \in Im\alpha$ if and only if $(a_{7})$ has an integer solution.

\item $p \in Im\overline{\alpha}$ if and only if $(b_{2})$ has an integer solution. In this case,

$p\equiv 1\pmod8, (\frac{p}{q})=+1$.

\item $q \in Im\overline{\alpha}$ if and only if $(b_{3})$ has an integer solution. In this case,

$q\equiv 1\pmod8, (\frac{q}{p})=+1$.

\item $pq \in Im\overline{\alpha}$ if and only if $(b_{6})$ has an integer solution. In this case,

$p, q \equiv1\pmod8$.

\item $(b_{1}), (b_{4}), (b_{5})$ and $(b_{7})$ does not have any integer solution.

\item If $p, q\equiv 5 \pmod8$ then $(a_{2}), (a_{3}), (a_{4})$ and $(a_{5})$ does not have any integer solution.
\end{enumerate}
\end{cor}
\begin{cor}
By the previous corollary, we have:
\begin{enumerate}
\item Let $p\equiv1, q\equiv5\pmod8$ and $(\frac{p}{q})=-1$.

$(b_{i}), 1\leq i\leq7$, does not have any integer solution. Hence $Im\overline{\alpha}= \lbrace1\rbrace$.

By the quadratic reciprocity, we have $(\frac{q}{p})=-1$.

$(\frac{-p}{q})=(\frac{p}{q})=-1$ and $(\frac{-q}{p})=(\frac{q}{p})=-1$.

It is possible for $(a_{7})$ to have integer solutions. So $Im\alpha \subseteq \lbrace \pm1, \pm 2pq \rbrace$.

As \eqref{eq2}, implies $\mid Im\alpha\mid\geq 4$. Consequently $Im\alpha=\lbrace\pm 1, \pm2pq\rbrace$.

Hence $r = 0$, i.e., $2pq$ is not a congruent number.
\\
\item Let $p\equiv 1, ~q\equiv3, 7\pmod8$ and $(\frac{p}{q})=-1$.

$(b_{i}), 1\leq i\leq7$, does not have any integer solutions. Hence $Im\overline{\alpha}= \lbrace 1\rbrace$.

By the quadratic reciprocity, we have $(\frac{q}{p})=-1$.

$(\frac{-q}{p})=(\frac{q}{p})=-1$ and $(\frac{\pm2q}{p})=(\frac{\pm2}{p}) \cdot (\frac{q}{p})=-1$.

It is possible for $(a_{1}), (a_{6})$ and $(a_{7})$ to have integer solutions. So $Im\alpha\subseteq \lbrace \pm1, \pm2, \pm pq, \pm2pq \rbrace$.

Since $2pq $ is a congruent number, (see\cite {mon}), then $r\geq1$.

As \eqref{eq2}, implies $Im\alpha= \lbrace \pm1, \pm 2, \pm pq, \pm 2pq \rbrace$. Hence $r=1$.
\\
\item Let $p\equiv 5\pmod8$ and $q\equiv3, 7\pmod8$.

$(b_{i}), 1\leq i\leq7$, does not have any integer solution. Hence $Im\overline{\alpha}= \lbrace 1\rbrace$.

Since $2pq$ is a congruent number, (see\cite {mon}), then $r\geq 1$.

It is possible for $(a_{2}), (a_{3}), (a_{4}), (a_{5})$ and $(a_{7})$ to have integer solutions. So $Im\alpha\subseteq \lbrace \pm 1, \pm p, \pm 2q, \pm 2pq \rbrace$ or $ \lbrace \pm1, \pm q,\pm 2p, \pm 2pq \rbrace$.

 As \eqref{eq2}, implies $Im\alpha= \lbrace \pm 1, \pm p, \pm 2q, \pm 2pq \rbrace$ or $\lbrace \pm 1, \pm q, \pm 2p, \pm 2pq \rbrace$. Hence $r=1$.
\\
\item Let $p, q\equiv 5\pmod8$.

$(b_{i}), 1\leq i\leq7$, does not have any integer solution. Hence $Im\overline{\alpha}= \lbrace 1\rbrace$.

It is possible for $(a_{7})$ to have integer solution. So $Im\alpha\subseteq \lbrace \pm1, \pm2pq\rbrace$.

As \eqref{eq2}, implies $Im\alpha= \lbrace \pm1, \pm2pq\rbrace$.

Hence $r= 0$, i.e., $2pq$ is not a congruent number.
\end{enumerate}
\end{cor}


\begin{thebibliography}{HD}
\bibitem{B1}
B. J. Birch, Diophantine analysis and modular functions, International Colloquium on Algebraic
Geometry, Tara Institute Studies in Mathematics 4, 3542 (1968).
\bibitem{coh}
H. Cohen, Number Theory. Vol. I: Tools and Diophantine Equations, Graduate Texts in Mathematics, 239 (Springer, New York, 2007).
\bibitem{hee}
K. Heegner, Diophantische analysis und modulfunktionen. Math. Z. 56, 227-253 (1952).
\bibitem{mon}
P. Monsky, Mock Heegner Points and Congruent Numbers, Math. Z. 204, 45 68 (1990) Mathematische
Zeitschrift, 9 Springer-Verlag 1990.
\bibitem{nag}
T. Nagell, L’analyse indetermineede degre superieur,´ Gauthier-Villars, Paris,(1929), 39,16 -17.
\bibitem{sil}
A. Silverman, Open questions in Arithmetic Geometry (Park, City ut,1999), IAS/Par Mathamatics Series q, AMS, Providence, RI (2001), 85-142.
\bibitem{ste}
N. M. Stephens, Congruence properties of congruent numbers. Bull. Lond. Math. Soc. 7, 182-184
(1975)
\bibitem{tun}
J. B. Tunnell, A classical diophantine problem and modular forms. Invent. Math. 72, 323-334
(1983).
\bibitem{was}
L. c. Washington, Elliptic Curves: Number Theory and Cryptography, Chapman- Hall, 2008,
\bibitem{Adl}
A. Adler, J. E. Coury, The Theory of Numbers : a Text and Source Book of Problems, Jones and Bartlett Publishers,1995.
\end{thebibliography}
\end{document}